\documentclass[a4paper,10pt,preprint]{elsarticle}
\usepackage[utf8]{inputenc}
\usepackage[english]{babel}

\usepackage{amssymb}\usepackage{amsmath}
\usepackage{graphics}
\usepackage{color}
\usepackage{bm}
\usepackage{url}
\newtheorem{thm}{Theorem}
\def\capacity{\mathop{\rm cap}\nolimits }
\def\Re{\mathop{\rm Re}\nolimits }
\def\Im{\mathop{\rm Im}\nolimits }
\def\be{\begin{equation}}
\def\ee{\end{equation}}
\def\bis{\hskip.5em\relax }

\begin{document}
%\headers{Capacity of symmetric polygonal condensers}{S. Bezrodnykh, A. Bogatyr\"ev, et al.}
\title{On capacity computation for symmetric polygonal condensers}
\author[CCAS,RUDN]{Sergei Bezrodnykh}
\ead{sbezrodnykh@mail.ru}

\author[INM]{Andrei Bogatyr\"ev\corref{cor1}}
\ead{ab.bogatyrev@gmail.com}

\author[INM]{Sergei Goreinov}
\ead{sergei.goreinov@ya.ru}

\author[INM]{Oleg Grigoriev}
\ead{guelpho@mail.ru}

\author[AA]{Harri Hakula}
\ead{harri.hakula@aalto.fi}

\author[TU]{Matti Vuorinen\corref{cor1}}
\ead{vuorinen@utu.fi}

\address[INM]{Institute for Numerical Mathematics,
Russian Academy of Sciences, Moscow, Russia}
\address[CCAS]{Computer Center of Russian Academy of Sciences, Moscow, Russia}
\address[AA]{Aalto University, Department of Mathematics and Systems Analysis, Espoo, Finland}
\address[TU]{Department of Mathematics and Statistics, University of Turku, Turku, Finland}
\address[RUDN]{RUDN University, Moscow, Russia}

\cortext[cor1]{Corresponding author}
\date{}

\begin{abstract}
Making use of two different analytical-numerical methods for capacity computation, we obtain 
matching to a very high precision numerical values for capacities of a wide family of planar condensers. These two methods are based respectively on the use of the Lauricella function \cite{bevl1,be1,be2} and Riemann theta functions \cite{B, gr, BG1, BG2}. We apply these results 
to benchmark the performance of numerical algorithms, which are based on adaptive 
$hp$--finite element method \cite{hrv1, hrv2, hrv3}  and boundary integral method 
\cite{ts, js, ac}. 
\end{abstract}


\maketitle
\section{Introduction}
Given a multiply connected domain $\mathcal A$ in the complex plane (more generally,
on a Riemann sphere) with piecewise smooth boundary, consider
a function $u(z)$ 
 harmonic in $\mathcal A$ 
with locally constant boundary values: $u=0$ at the outer boundary component and $u=1$
at all inner boundary components. The Dirichlet integral of this function over $\mathcal A$ is
exactly the capacity of the `condenser' $\mathcal A$
(see e.g. \cite{ah, avv, abr, js, ts}). 
The study of capacity attracts much attention due to the place it occupies in many 
theoretical and applied questions. Thus, the capacity or the close concept of a conformal 
modulus plays an important role in several issues of function theory, its analytic, 
geometric problems and problems of conformal and quasiconformal mappings 
(see e.g. \cite{ab,avv,du,ku,posz}).

Besides, condensers or conformal moduli occur in several problems of 
physics, mechanics and engineering. A general description of such applications is 
given in the monographs \cite{he, ku, ky, ps, sl}. The traditional applications
include the calculation of electrical resistance, capacity, and other 
physical characteristics of polygonal--shaped elements of integrated circuit networks 
(see e.g. \cite{da, ga95, ll, mi, p89, rc, tic, tr, thg}). 
In recent decades, a number of unconventional applications 
have emerged, including simulation in electrical impedance tomography, airfoil modeling 
in computational fluid dynamics, in computer vision and other issues (see \cite{hyv, 
sm, ku, sl, thg} and the bibliography therein).

Despite the high demand for reliable results in this area, the development of relevant
effective and high-precision methods remains an actual problem far from the exhaustive 
solution. The papers \cite{bsv,hrv1,hrv2,hrv3} have extensive bibliographies which more 
or less cover what is known about numerical computation of capacities of condensers 
or moduli of quadrilaterals. However, thorough observation of this area shows that 
not many numerical results are known. Perhaps the longest list of such results is given 
in the paper \cite{bsv} and in the books \cite{ps, ku, sl} where references to earlier literature 
may be found.
At present, there are two groups of numerical methods for treatment of this problem: 
direct and indirect ones. Direct ones perform the capacity computation by solving
a boundary value problem or by numerical implementation of a relevant variational 
principle. Indirect methods transplant our problem by conformal mapping onto 
a new domain where the problem is trivial. So those methods are in fact the methods 
of conformal mapping. 

One of the most popular direct method is the finite element method (FEM). It was 
first applied to conformal modulus computation in \cite{o67,mi,ga72}, and then 
experienced a great development and generalization \cite{w79,p89,ham,ga95,thg}. 
In \cite{w79,bo}, it had been improved by implementation of `singular functions'.
In \cite{bsv} an adaptive finite element method (AFEM) was proposed. General 
progress in the development of FEM was represented in 
monographs \cite{szba} and \cite{sw}. The implementation of these approaches to 
the computation of capacities was carried out in the works \cite{hrv1, hrv2, hrv3}, 
where $hp$--FEM was improved. In the present work this method, along with other methods, 
is used for obtaining numerical results. 

It is known that FEM gives an upper bound for the Dirichlet integral, and consequently 
for the capacity.
The Gauss--Thomson extremal principle \cite{hps, ts, kir, w79, ga95}, dual to the Dirichlet principle underlying 
the finite element method, can also be used to provide lower bounds for the capacity.
In \cite{w79}, the `singular functions' technique was 
introduced in FEM. 
A discrete analogue of this method was developed \cite{ts, kir},
which uses point systems similar to the Fekete system arising in connection 
with the notion of transfinite diameter of a compact set.
In works \cite{sk, p89}, the finite difference method has been used to calculate 
the capacities. In doing so, the Bakhvalov--Fedorenko multigrid and other special techniques were applied 
to improve its efficiency. 
Boundary integral equation methods have been also applied for capacity computation 
\cite{js, ac, rei, p89, cwlh}, with the use of `singular functions'. In the present work this method, 
along with others, is used for obtaining numerical results.
Another direct method of capacity computation without determining the conformal map 
itself was developed in \cite{ga64, lps, pwh, ppss, mss, p03}. The method is based on 
the Bergman kernel, the Bieberbach polynomials and orthonormalization procedure 
with the use of `singular functions' technique. Finally, for a special class of 
elongated domains, the domain decomposition method was developed \cite{hot, 
gah, l}.

Turning to indirect methods, which are based on approximating 
the conformal map, we must mention monographs and surveys \cite{ga64,tr86,p89,rssv,he,
ky,dt,p03,sl,ku,ps}. This group of methods includes Symm's method developed in 
\cite{js,ac,h,p89,hlcw,cwlh,rssv,p03} and the Theodorsen--Garrick method and 
its modification such  as the Fourier series method (see e.g. 
\cite{tg,ga64,p89,rssv, ky,he}). For polygonal domains, undoubtedly, the main method is 
the Schwarz--Christoffel (SC) integral (see e.g. \cite{ga64,tr80,tr86,rssv,he,ky,he,sl}) 
and its generalizations to multiply connected domains \cite{da,cr,dek,deek,et}. Of particular 
note is Driscoll--Trefethen SC Toolbox \cite{dr} based on compound Gauss--Jacobi 
quadrature \cite{tr86,dt} and Hu's software for doubly connected domains \cite{hu}.

Recent results on the analytic continuation of the Lauricella function \cite{be1,be2,be3} and 
applications of
the theta functions theory \cite{B} open new opportunities in the development of SC--technology 
(\cite{bevl1,be3} and \cite{gr,BG1}, respectively). In the present paper, the efficiency of the two 
corresponding methods is demonstrated in the calculation of the capacity of several families of 
elaborate polygonal domains. We apply these results to benchmark the performance of numerical 
algorithms, based on adaptive $hp$--finite element method \cite{hrv1, hrv2, hrv3} and boundary 
integral method \cite{ts, js, ac}.

{\bf Acknowledgements.} The work of AB, SG and OG was supported by RSCF grant 16-01-10349, SB acknowledges 
the receipt of support from the RUDN University 'Program 5-100' and RFBR grant 19-07-00750.
All authors are indebted to Professors Vladimir Vlasov and Antti Rasila who supported this research at all its stages.

\section{Problem setting}

\label{section2}

We consider planar polygonal condensers in the complex plane of variable $z = x + i y$
such that the calculation of their capacity can be carried out with the aid of the
SC integral. Let $\cal A$ be a polygonal domain of finite connectivity in the plane. 
We designate the outer boundary of polygon $\cal A$ as $\partial_0\cal A$ and the collection of all the rest  
boundaries as $\partial_1\cal A$. The capacity of condenser $\cal A$ is the Dirichlet (or energy) integral 
\be
\label{CapDef}
\capacity({\cal A})=\int_{\cal A} |\nabla ~U(z)|^2~dxdy
\ee
calculated at the solution  $U (z)$ of the  Dirichlet boundary value problem:  
\begin{equation} \label{DirProbl}
\Delta U (z) = 0,\quad z\in{\cal A};
\qquad U (z) = 0,\quad z \in\partial_0 {\cal A};
\qquad U (z) = 1,\quad z \in\partial_1 {\cal A}.
\end{equation}
This solution $U (z)$ is understood in the classical sense, namely $U\in C (\overline{\cal A}) \cap C^2 (\cal A)$.

\begin{figure}[t]
%\vspace*{-1.5cm}
 \begin{center}
 \begin{picture}(35, 30)
\thinlines
\multiput(-5,10)(2,0){30}{\line(1,0){1}}

\thicklines
\put(5,25){\line(1,0){30}}
\put(5,10){\line(0,1){15}}
\put(35,10){\line(0,1){15}}

\put(10,10){\line(0,1){10}}
\put(10,20){\line(1,0){15}}
\put(25,15){\line(0,1){5}}
\put(25,15){\line(1,0){5}}
\put(30,10){\line(0,1){5}}

\put(31,5){$N$}
\put(6,5){$N$}
\put(19,16){$F_1$}
\put(10,26){$F_0$}
\put(29,19){$\cal D$}
\put(40,15){$\mathbb{H}$}

\end{picture}
\caption{Half the condenser $\cal A$.}
\label{HalfCond}
\end{center}
\end{figure}
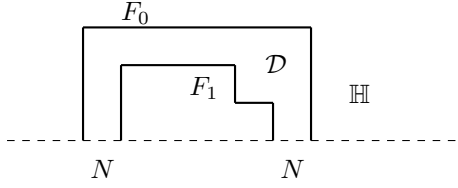

In what follows we assume that  $\cal A$ is mirror symmetric with respect to the real axis and its intersection with
the upper half-plane $\mathbb{H}\, := \,\{ z: \Im z > 0\}$ is a simply connected 
domain ${\cal D}\, :=\, {\cal A}\, \cap \mathbb{H}$.  Its boundary
$\partial{\cal D}$ is the union $\partial {\cal D} = F_0 \cup F_1 \cup N$ of
compact sets $F_s:=\partial_s {\cal A} \cap ({\mathbb H}\cup{\mathbb R})$, $s=0,1$, 
and $N := {\cal A}\cap {\mathbb R}$ -- see Fig.~\ref{HalfCond}. From the symmetry of the domain $\cal A$ it follows that the Dirichlet problem (\ref{DirProbl})
can be reduced to the following mixed boundary value problem in the domain $\cal D$:
\begin{equation} \label{DirProb2}
%\begin{array}{rl}
\Delta U (z) = 0,  \quad z \in {\mathcal D};
\qquad U (z) = 0 ~\hbox{(resp. $1$)}, \quad z \in F_0~\hbox{(resp. $F_1$)};
%\qquad u (z) = 1, \quad z \in F_1;
\qquad\partial_n U (z) =0, \quad z \in N,
%\end{array}
\end{equation}
where $\partial_n$ denotes the outer normal derivative; solution of the reduced problem
(\ref{DirProb2}) is also understood in the classical sense. So, the capacity
(\ref{CapDef}) can be evaluated via the solution to the problem (\ref{DirProb2})
by the formula
\begin{equation} \label{capred1}
\capacity({\cal A}) = 2\int_{\cal D} |\nabla u(z)|^2 ~dxdy.
\end{equation}

\section{Description of four approaches to capacity computation}  

\subsection{Capacity of a rectangular condenser and theta functions}
\label{Bogatyrev}
In this section we briefly describe how to get the effective closed formulas for the capacity
of an axisymmetric condenser bounded by a rectangular chain. The approach is based on two observations.
First, the Dirichlet integral (or the energy of the electric field inside the capacitor) is invariant under conformal mappings.
Second, the SC integral that conformally maps the upper half plane to a simply connected rectangular polygon
is an abelian integral on a suitable hyperelliptic Riemann surface. Here for simplicity we consider the case when the surface is elliptic, more instructions on higher genus curves may be found in \cite{B,BG1,BG2,gr}.

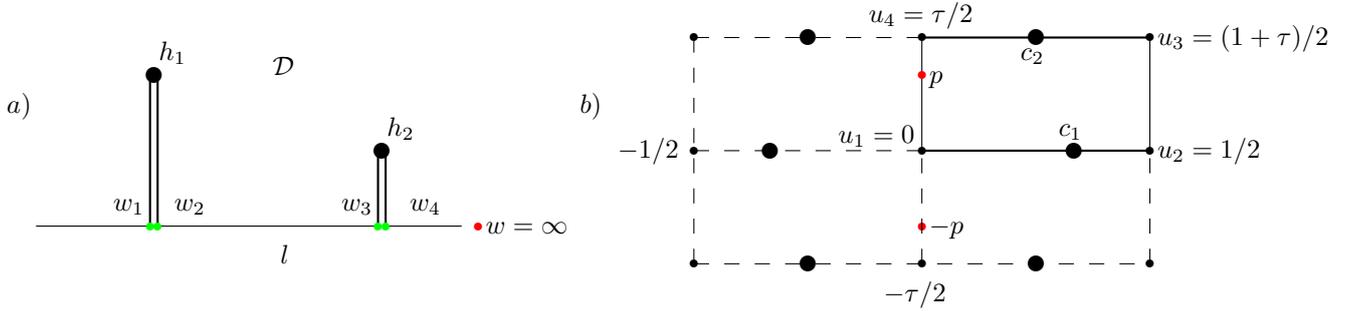
\begin{figure}[t]
\vspace*{2ex}
\begin{picture}(60, 30)
\linethickness{.5mm}
\thinlines
\put(5,5){\line(1,0){15}}
\put(21,5){\line(1,0){29}}
\put(51,5){\line(1,0){10}}

\thicklines
\put(20,5){\line(0,1){20}}
\put(21,5){\line(0,1){20}}
\put(50,5){\line(0,1){10}}
\put(51,5){\line(0,1){10}}
\put(20.5,25){\circle*{2}}
\put(50.5,15){\circle*{2}}
{\color{green}
\put(20,5){\circle*{1}}
\put(21,5){\circle*{1}}
\put(50,5){\circle*{1}}
\put(51,5){\circle*{1}}}

{\color{red}
\put(62,5){\circle*{1}}
}

\put(63,4){$w=\infty$}
\put(14,7){$w_1$}
\put(22,7){$w_2$}
\put(44,7){$w_3$}
\put(53,7){$w_4$}
\put(0,20){$a)$}
\put(36,0){$l$}
\put(20,27){$h_1$}
\put(50,17){$h_2$}
\put(35,25){$\cal D$}
\end{picture}
\hspace{2.5cm}
\begin{picture}(100, 30)
\linethickness{.5mm}
\thicklines
\put(35,15){\line(1,0){30}}
\put(35,30){\line(1,0){30}}
\thinlines
\put(65,15){\line(0,1){15}}
\put(35,15){\line(0,1){15}}

\put(55,15){\circle*{2}}
\put(15,15){\circle*{2}}

\put(53,17){$c_1$}
\put(48,27){$c_2$}

\put(50,30){\circle*{2}}
\put(50,0){\circle*{2}}

\put(20,30){\circle*{2}}
\put(20,0){\circle*{2}}

{\color{red}
\put(35,25){\circle*{1}}
\put(35,5){\circle*{1}}
}

\put(36,24){$p$}
\put(36,4){$-p$}

\multiput(5,0)(30,0){3}{\circle*{1}}
\multiput(5,15)(30,0){3}{\circle*{1}}
\multiput(5,30)(30,0){3}{\circle*{1}}

\multiput(5,0)(4,0){15}{\line(1,0){2}}
\multiput(5,15)(4,0){8}{\line(1,0){2}}
\multiput(5,30)(4,0){8}{\line(1,0){2}}
\multiput(5,0)(0,4){8}{\line(0,1){2}}
\multiput(35,0)(0,4){4}{\line(0,1){2}}
\multiput(65,0)(0,4){4}{\line(0,1){2}}

\put(24,16){$u_1=0$}
\put(66,14){$u_2=1/2$}
\put(66,29){$u_3=(1+\tau)/2$}
\put(28,32){$u_4=\tau/2$}
\put(-5,14){$-1/2$}
\put(30,-5){$-\tau/2$}

\put(-10,20){$b)$}
\end{picture}

\vspace{3mm}
\caption{Half the condenser $\cal D$ (left) is mapped to the rectangle $\frac12\times\frac{|\tau|}2$ (right).}
\label{F2Rect}
\end{figure}

As a simple illustration we consider here the case of two parallel slots of lengths $2h_1$, $2h_2$ at distance $l$ with the common symmetry axis perpendicular to the slots -- see condenser $E$ in Fig.~\ref{DEFG}. The symmetry axis divides the condenser into two halves as shown on the left  Fig.~\ref{F2Rect}.  The upper half $\cal D$ of the capacitor may be conformally mapped to the upper half plane and the image of the four vertexes $w_1,\dots,w_4$ at its right angles give us four branch points of an elliptic surface where the SC integral delivering the inverse mapping is properly defined.

However, we use another model of this surface, namely the factor of the complex plane of the variable $u$ by the lattice $\mathbb{Z}+\tau\mathbb{Z}$ with $\tau\in i\mathbb{R}^+$. The absolute value $|\tau|$ of the elliptic modulus appears as the aspect ratio of the rectangle conformally equivalent to $\cal D$ with the banks of the slots  being mapped to the horizontal sides of the rectangle. In this model of the condenser, the equilibrium field \eqref{DirProb2} inside it is proportional to $\Im u$ and we immediately get that 
\be
\label{CapSlots}
\capacity({\cal A})=2/|\tau|,
\ee
$|\tau|$ being the ratio of the height of the rectangle to its width. Now we derive a system of transcendental equations 
for the elliptic modulus $\tau$.

We want to get an expression for the SC integral $w(u)$ in terms of the variable $u$ of the universal cover of the torus. We start with the building block for the function theory on a torus (and other closed surfaces) -- the third kind abelian integral $\eta_{pq}$ that has logarithmic poles at two given points $p,q$ of the surface:
$$
\eta_{pq}(u)=\log\frac{\theta_1(u-p)}{\theta_1(u-q)},
\qquad 
\theta_1(u):=\theta_1(u|\tau):=2\sum_{k=0}^\infty(-1)^k\exp(\pi i(k+1/2)^2\tau)\sin(\pi (2k+1)u),
$$
$\theta_1(u)$  being the only odd elliptic theta function of modulus $\tau$. Theta functions are represented by extremely quickly convergent series. Now we use the transformation properties of theta functions \cite{A, Mum} and check that this integral is normalized with respect to the horizontal meridian of the torus ($a-$ cycle) and has a nontrivial period along the vertical meridian ($b-$ cycle):
$$
\eta_{pq}(u+1)-\eta_{pq}(u)=0;\qquad
\eta_{pq}(u+\tau)-\eta_{pq}(u)=2\pi i(p-q).
$$
One may construct an $a-$ normalized second kind integral with first order pole at the point $p$ as the only singularity on the torus just differentiating the above integral $\eta_{pq}$ with respect to the first pole position:
$$
\omega_{2p}:=-\frac{\theta'_1(u-p)}{\theta_1(u-p)}=-1/(u-p)+locally~holomorphic~term
$$

From an explicit algebraic expression for SC integral $w$, one can check that the abelian differential $dw$ has exactly two double poles interchanged by the involution $Ju=-u$ of the surface fixing its four branch points $u=0,1/2, \tau/2, (1+\tau)/2$. Further analysis of geometry 
shows that $a-$ period of $dw$ vanishes while its $b-$ period is equal to $2l$.  Hence, SC integral takes the appearance
$$
w(u)=-\frac{l}{2\pi i}\left(
\frac{\theta'_1(u-p)}{\theta_1(u-p)}+\frac{\theta'_1(u+p)}{\theta_1(u+p)}
\right)
$$
with pole $p\in[0,\tau/2]$ not known at the moment. The system of equations for $p$ and $\tau$
includes also the positions of zeros $c_1\in[0,1/2]$ and  $c_2\in[0,1/2]+\tau/2$ of the differential $dw$
which correspond to the tips of the slots (fat dots in Fig.~\ref{F2Rect}b) ). Closed system of transcendental equations 
with respect to four real numbers $|\tau|>0$, $\Im p\in[0,|\tau|/2]$ and $c_1$, $\Re c_2\in[0,1/2]$
includes geometric parameters of the condenser and  consists of four (essentially) real equations:
\be
\label{SystemE}
dw(c_j)=0;\qquad
\Im w(c_j)=h_j;\qquad j=1,2.
\ee
By the use of the uniqueness of a normalized conformal mapping, one can show that this system has a unique solution in the above mentioned limits provided all the input numbers $l,h_1,h_2$ are positive. The system (\ref{SystemE}) may be solved by Newton method (with parametric continuation if necessary) -- see the discussion in \cite{B}. Summing up, to get the capacity of condenser $E$ of Fig.~\ref{DEFG},
one solves the system of four  transcendental equations (\ref{SystemE}) involving the geometrical dimensions of the condenser. One of the 
independent variables of the system, namely $|\tau|$ immediately gives the capacity via the formula (\ref{CapSlots}).

For more involved axisymmetric rectangular condensers the scheme of the capacity calculation differs only in technical details. The SC integral which maps the half-plane to half a condenser is a hyperelliptic one. It has the explicit expression in terms of Riemann theta functions. The latter depend on the variables in the Jacobian of the curve -- multidimensional complex torus defined by the period matrix $\tau$ of the associated curve. For the curve with genus $g>1$  the Jacobian is different from the curve itself and we need to characterize the locus of the curve inside its Jacobian \cite{B2}. For $g>2$ we have to characterize the locus of the hyperelliptic period matrices in the Siegel space -- this is the one of the faces of the notorious Schottky problem \cite{Mum}.

\subsection{Lauricella function Approach to the Solution of Problem}
\label{Vlasov}
This method of capacity computation is based on using the SC integral technique
\cite{bevl1,bevl2} in combination with new results on the Lauricella
function theory \cite{be1, be2,be3}. In the mentioned works the method was successfully
applied to conformal mapping of a number of complex--shaped domains with highly efficient
solution of parameter problem for the SC integral in crowding situation
(very close location of pre-vertices).

For the sake of simplicity we describe the method in the case when the condenser $\cal A$ 
is a topological annulus. As in the previous section we can conformally map 
half the condenser $\cal D$ to a rectangle with the boundaries $F,N$ being mapped to the 
vertical/horizontal sides respectively. The mixed boundary value problem \eqref{DirProb2}
now becomes trivial and  we arrive to formula \eqref{CapSlots} for the capacity of $\cal A$
with the aspect ratio $|\tau|$ being the ratio of width of the rectangle to its height. Generalization of this
construction for higher connectivity condensers may be found in \cite{ks}.

The first step toward finding of the value of $|\tau|$ is the construction of a conformal mapping 
of the upper half-plane $\mathbb{H}$ onto the domain $\mathcal{D}$. 
We denote the vertices of polygonal domain $\,\mathcal{D}\,$ by $\,z_k$, $\,k = 1, 2, \ldots, K,\,$ 
and values of angles in these vertices by $\pi\,\alpha_k$. We denote the pre--vertices of $z_k$ under this mapping by 
$\zeta_k: = \Phi^{-1} (z_k)$ and suppose, that $\zeta_1 = \infty$, $\zeta_2 = 0$, $\zeta_3 = 1$.
In other words, we subject function $z = \Phi (\zeta)$ to the conditions
\begin{equation} \label{condForPhi}
\Phi (\infty) = z_1\,, \qquad \Phi (0) = z_2\,, \qquad \Phi (1) = z_3\,,
\end{equation}
that ensure its existence and uniqueness. So, mapping $\Phi (\zeta)$ can be presented 
in the form of the following SC integral:
\begin{equation} \label{SC-FormIni}
\Phi (\zeta) = C_0 \int_0^\zeta \, \prod_{j = 2}^K \,(t - \zeta_j)^{\alpha_j -1} d t\, 
+\,z_2\,,
\end{equation}
where $\zeta_2 = 0$, $\zeta_3 = 1$. Here the factor in the integrand corresponding to $\zeta_1$  
is omitted because  $\zeta_1 = \infty$. 

Presentation (\ref{SC-FormIni}) contains $(K - 2)$ unknown parameters, namely, the prevertices 
$\zeta_4, \zeta_5,\ldots, \zeta_K$ and the coefficient $C_0$. For finding these parameters,
we compose  in conventional way \cite{ks} the following system of $(K - 2)$
transcendental equations with $(K - 2)$ unknowns:
\begin{equation} \label{SyTrascEq}
\Big| C_0 \int_{\zeta_k}^{\zeta_{k + 1}} \, \prod_{j = 2}^K \,(t - \zeta_j)^{\alpha_j - 1} d t\,\Big| =
\, |z_{k + 1} - z_k|\,,\quad k = 2, 3,\ldots , K - 1\,.
\end{equation}
The solution of this system is carried out, as usually, by means of the Newton method. 
Since this method is iterative, then for providing its efficiency the following 
two conditions should be fulfilled: 1) high precision computation of integrals 
in expression (\ref{SyTrascEq}) ought to be performed; 2) a good  initial approximation 
must be constructed.

For computing the integrals of interest we make use of some new results \cite{be1, be2, be3}
from the theory of the Lauricella function $F_D^{(N)}$, which is a hypergeometric 
function of $N$-dimensional complex variable $\mathbf{x} := (x_1, \ldots, x_N)$. 
The basic point is the Euler-type integral representation \cite{e} for this function:
\begin{equation}\label{LF_2a}
F_D^{(N)}\, (\mathbf{a}; b, \, c;\, \mathbf{x})\,
=\,
\frac{\Gamma (c)}{\,\,\Gamma(b)\, \Gamma (c - b)}\,
\int\limits_0^1\,
\frac{\,\,t^{b - 1}\, (1 - t)^{c - b - 1}}{\prod\limits_{j = 1}^N\,
(1 - t\, x_j)^{\,a_j}}\,\, dt\,,
\end{equation}
where $\mathbf{a} := (a_1, \ldots, a_N)$, $b$, $c$ are parameters, $\Gamma(t)$ the 
gamma function. It can be easily noticed that by the use of representation (\ref{LF_2a}) 
the integrals in (\ref{SyTrascEq}) can be rewritten in terms of the Lauricella function. 
In recent works \cite{be1, be2, be3}, the formulas of analytic continuation for 
the Lauricella function have been derived; for any point $\mathbf{x} \in \mathbb{C}^N$, 
these formulas give appropriate exponentially convergent series for 
$F_D^{(N)}\, (\mathbf{a}; b, \, c;\, \mathbf{x})$. Thus, we obtain high precision 
algorithm of computation of integrals in expression (\ref{SyTrascEq}). So, the first 
condition of efficiency of Newton's method is satisfied.

For a number of elaborate domains \cite{bevl1,bevl2}, 
it was found a good initial approximation for system (\ref{SyTrascEq}) of 
transcendental equation for parameters of the SC integral. 
It was made by application of the theory of conformal mapping of singularly 
deformed domains. 
The base for such application lies in the fact, that a complex-shaped domain very 
frequently can be considered as a singularly deformed domain of much simpler
configuration, so that for the latter domain a conformal mapping can be easily constructed.
The so called `deformation parameter' also can be pointed out, and in this case, 
as a rule, its numerical value is small. The mentioned theory gives an asymptotics 
for the conformal mapping as the deformation parameter goes to zero. This asymptotics
gives the required initial approximation.

Applying the described algorithm for integral computation in combination with obtained 
initial approximation ensures high accuracy computation of unknown parameters in
representation (\ref{SC-FormIni}). When implementing the Newton iteration process, it is necessary that at every step the approximate values $\zeta_j$, $j = 4, \dots, K$, additionally satisfy the conditions $\zeta_{j + 1}> z_j$, $j = 3, \dots, K-1$. This difficulty is indicated in \cite{tr80}, where a method for its elimination was proposed by transition to the so-called unconstrained variables. At the same time, if the initial approximations used for the prevertices are sufficiently close to their exact values, then the indicated order conditions are satisfied. The asymptotics we used allows to construct good initial approximation, which ensures that the order is preserved for the approximate values of the prevertices at each step of the iterative algorithm. Thus, the conformal mapping $\,z = \Phi (\zeta)\,$ is 
constructed.  Next we find the preimages $\zeta^{(j)}$ under this mapping of four endpoints 
$z^{(j)}$ of two segments $N$ (see Fig \ref{HalfCond}) enumerated in increasing order starting from the smallest 
$z^{(1)}:=F_0\cap N_0$.

Now we transform by a M\"obius mapping the 4-tuple $(\zeta^{(1)}, \dots, \zeta^{(4)})$ to a standard form 
$(-1, 1, 1/k, -1/k)$ with a certain elliptic modulus $k\in(0,1)$. 
It is not difficult to verify that the mapping is given by the formula
\begin{equation} 
\frac{T (\zeta) + \sqrt{\varkappa}}{ - T (\zeta) + \sqrt{\varkappa}}\,,
\end{equation}
where
\begin{equation} \label{MobiusMap}
T (\zeta) = \frac{\zeta - \zeta^{(2)}}{\zeta - \zeta^{(1)}}\,\cdot\, 
\frac{\zeta^{(3)} - \zeta^{(1)}}{\zeta^{(3)} - \zeta^{(2)}}\,,\qquad
\varkappa = \frac{ \zeta^{(4)} - \zeta^{(2)}}{ \zeta^{(4)} - \zeta^{(1)}}\,\cdot\, 
\frac{\zeta^{(3)} - \zeta^{(1)}}{\zeta^{(3)} - \zeta^{(2)}}\,,
\end{equation}
and the parameter $k$ is given by the expression
\begin{equation}\label{ValueOfTau}
k = \frac{\sqrt{\varkappa} - 1}{\sqrt{\varkappa} + 1}\,.
\end{equation}

The first kind incomplete elliptic integral of the modulus $k\in(0,1)$ maps the upper halfplane to a rectangle, 
with the set $(-1, 1, 1/k, -1/k)$ being mapped to its vertexes \cite{A}.  The aspect ratio $|\tau|$ of this rectangle is twice the 
ratio of complete elliptic integrals $K(k),~K'(k)$ of the same modulus. Expressing complete elliptic integrals 
in terms of the Gauss hypergeometric function  $F (a, b; c; z)$  in the form $K(k) = (\pi/2) F(1/2,1/2;1;k^2)$
we eventually get the value of capacity
\begin{equation} \label{caprFin}
\capacity({\cal A}) = \frac{F \Bigl(\frac{1}{2},\, \frac{1}{2};\, 1; \,k^2\Bigr)}
{F \Bigl(\frac{1}{2},\, \frac{1}{2};\, 1; \,1 - k^2\Bigr)}\,.
\end{equation}

\subsection{hpFEM}
\label{sec:preliminaries}
In this section we briefly describe the high-order $p$-, and $hp$-finite element 
methods (FEM). The implementation used in this paper is based on our earlier work \cite{hrv1}. 
There exists an extensive literature on $hp$-FEM, our main references are the 
books by Szabo and Babu{\v{s}}ka, and Schwab, \cite{szba}, \cite{sw}, respectively.

In standard $h$-FEM the polynomial order of approximation $p$ 
is fixed and one refines the mesh with the purpose of improving accuracy.
For higher-order methods there are more options. It is also possible
to obtain a convergent sequence of solutions by keeping the mesh fixed but
varying the order of the elemental polynomial basis ($p$-version). 
The mesh may include triangular and quadrilateral
elements.
However, it is not necessary
to use the same polynomial basis over the whole domain. Indeed, the fundamental
idea behind the $hp$-version is to couple mesh generation and selection of the
desired elementwise polynomial order. 

For a certain class of problems it can be shown that if the mesh and the elemental
degrees have been set optimally, we can obtain \textit{exponential convergence}.
A geometric mesh is such that each successive layer of elements changes in
size with some \textit{geometric scaling factor} $\alpha$,
toward some point of interest. For planar problems the fundamental result
is given as follows:
\begin{thm} [{\cite{sw}}]\label{propermesh}
Let $\Omega \subset \mathbb{R}^2$ be a polygon, $v$ the FEM-solution, and
let the weak solution $u_0$ be in a suitable countably normed space where
the derivatives of arbitrarily high order are bounded.
Then
\[
\inf_v \|u_0 - v\|_{H^1(\Omega)} \leq C\,\exp(-b \sqrt[3]{N}),
\]
where $C$ and $b$ are independent of $N$, the number of degrees of freedom.
Here $v$ is computed on a proper geometric mesh, where the orders of individual
elements depend on their originating layer, such that the highest layers have the smallest
orders. The result also holds for constant polynomial degree distribution.
\end{thm}

The choice of both $h$ and $p$ will have an influence on the performance and 
the accuracy attained. The mutual influence of 
these choices is hard to analyze 
theoretically but usually it can be seen in the results.
In our experience the initial mesh or the first layer of elements
has to be designed with care in order to balance the errors
between the refined and non-refined parts corresponding to singularities
and the smooth part of the solution, respectively.
This process is explained in more detail in the following example.

\subsubsection{Example}

In the context of planar condensers the theory above is directly applicable.
The locations of the singularities are always known a priori.
We can always design the initial meshes so that 
geometric grading of the meshes can be carried out using a 
replacement rule-based approach. Notice, that standard meshing approaches
based on geometric data structures such as Delaunay triangulations cannot
support exponential grading. Furthermore, in our implementation we use rule-based meshing with exact arithmetic and hence
can generate meshes with geometric grading down to arbitrary small element sizes.
One should also notice that it is admissible to have cuts inside the domain or at 
the boundary. One should also notice that the replacement rules can lead
to meshes with different kinds of elements in successive levels of refinement.

\begin{figure}[t]
%\vspace*{-1cm}
  \centering
  \includegraphics[width=0.6\textwidth]{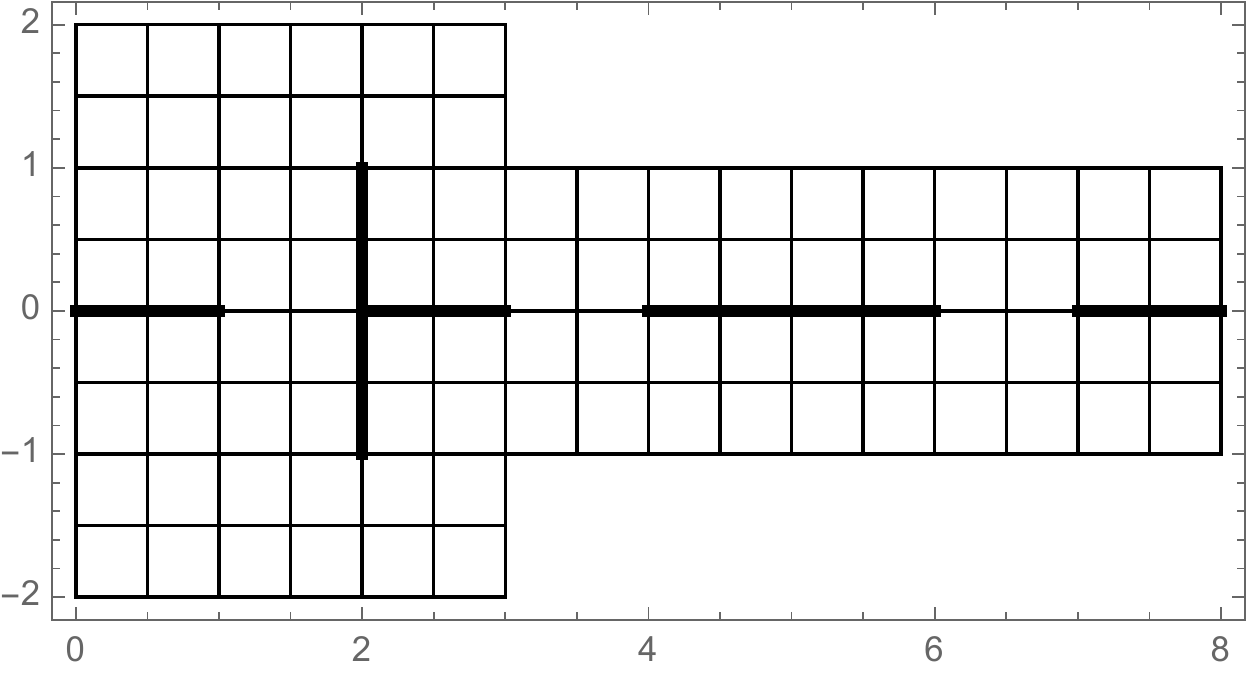} \\
  \includegraphics[width=0.6\textwidth]{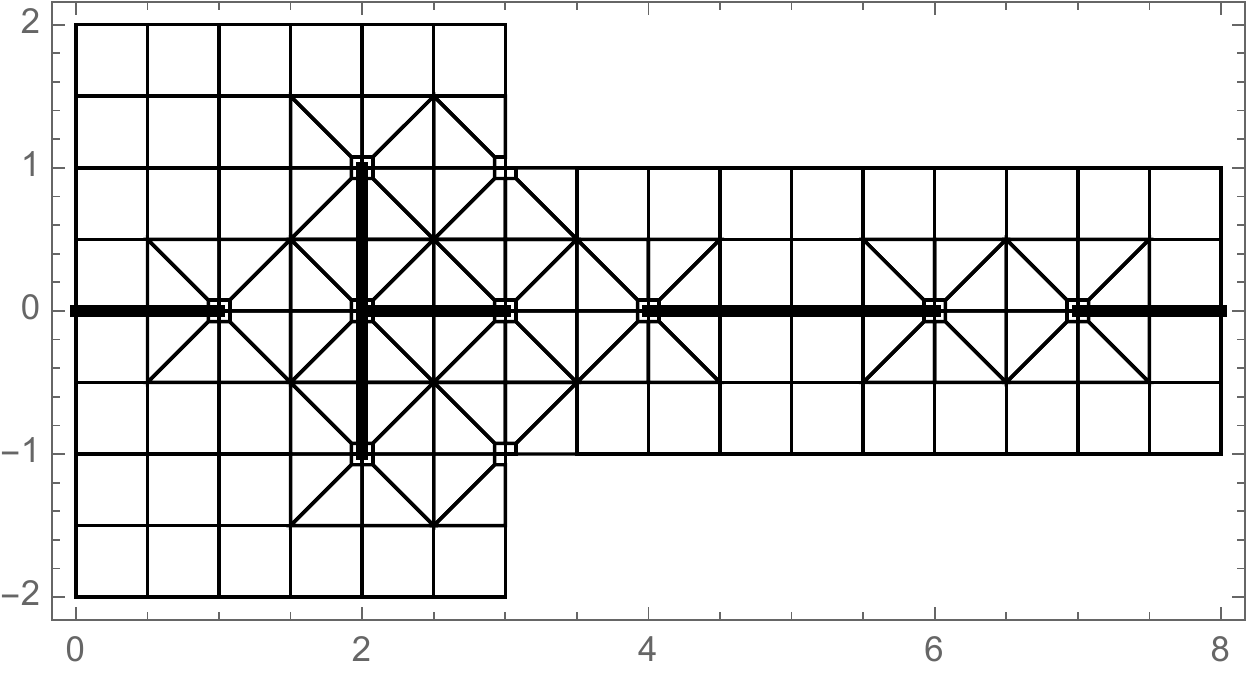} \quad
  \includegraphics[height=1.5in]{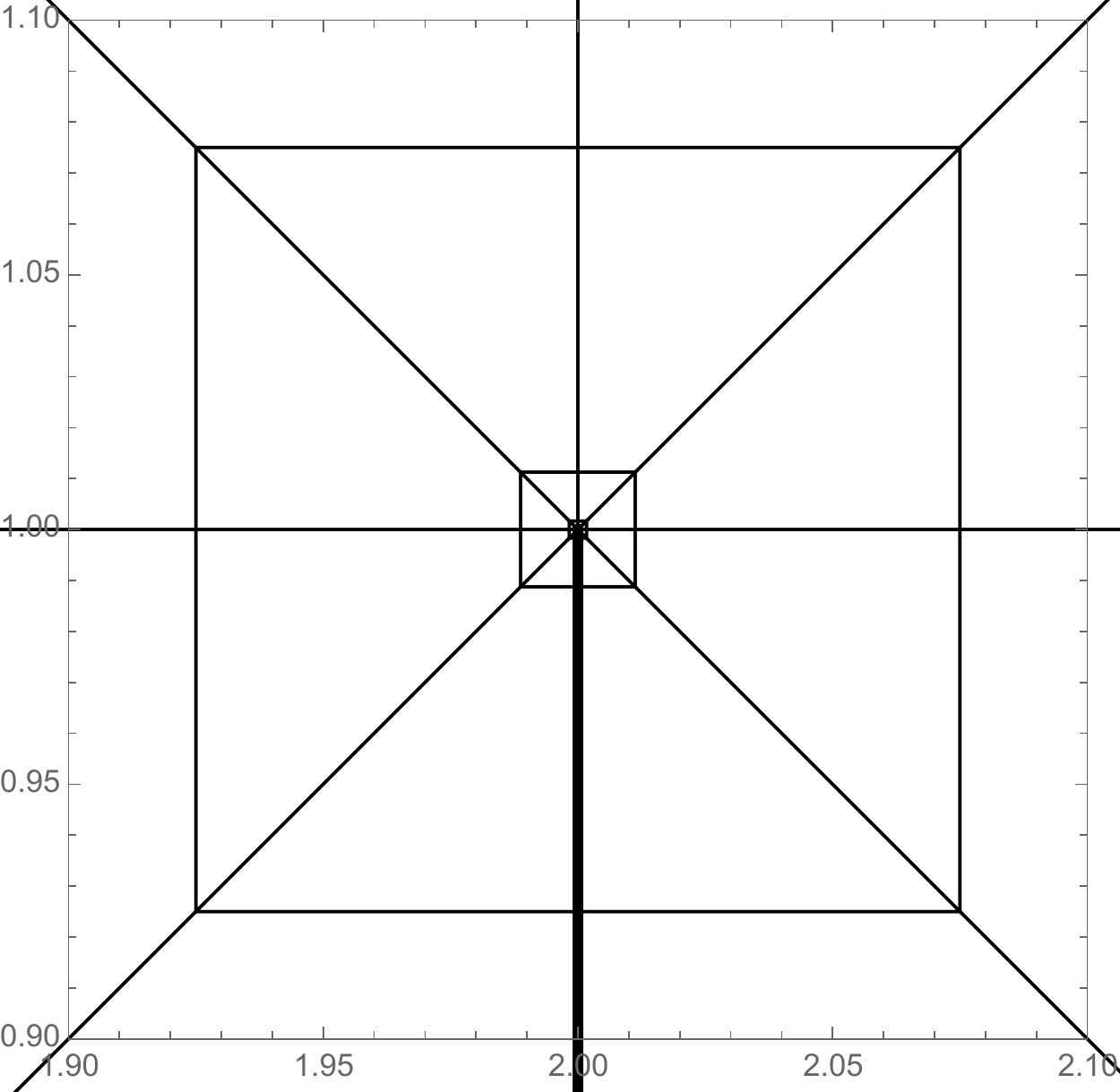}
  \caption{Mesh generation; From top counterclockwise: The initial mesh, 
  mesh after one level of refinement, detail of the mesh around the singularity at
  $(2,1)$ after 14 levels of refinement.}\label{fig:meshes}
\end{figure}

In Figure~\ref{fig:meshes} we show an example of the mesh generation
procedure. The initial mesh is laid out so that every singularity
and reentrant corner
is isolated, that is, the replacement rules are consistent over each
element adjacent to a singularity. This is evident in the mesh after
one level of refinement. If one level of refinement can be carried
out successfully, it follows immediately that the refinement process
can be repeated as many times as desired.

This kind of approach has a serious drawback, however. It is very difficult to
modify the meshes adaptively. Of course, changing the polynomial orders is straightforward.
Thus, 
instead of employing an adaptive algorithm we have simply modified the solution
scheme ($p$ and grading) based on a posteriori error estimation. 
In other words, the meshing process has been adapted rather than
any given mesh. The error estimation algorithm is described in \cite{Hakula2017}.

\subsection{Boundary integral equation approach}
Here we do not use the symmetry of ${\cal A}$. Denote the connected components of
$\partial_s{\cal A}$ as $\Gamma_j$,
$j=1, 2, \dots$,
and drop for a moment the explicit value
of the Dirichlet boundary condition. The whole boundary will be denoted
$\Gamma=\partial_0{\cal A}\cup\partial_1{\cal A}$.
Also, we disregard the subdomains
incident to only one particular boundary $\Gamma_j$, like inner/outer rectangles
and inner cross in (`swiss crosses') Fig.~\ref{DEFG}.

There are several ways to reduce the capacity computation to integral equations,
including single or double layer ansatz \cite{schwab}, or angular potential
leading to singular integral equations \cite{krut}.
If we define boundary integral operators $V$, $K$, $K'$, $W$ by
\begin{equation}\label{eq:layers}
\begin{aligned}
V\!\mu(z)&:=\int_\Gamma g(z-\zeta)\mu(\zeta)\,|d\zeta|,\quad 
K\!\mu(z):=\int_\Gamma \partial_{n(\zeta)}g(z-\zeta)\mu(\zeta)\,|d\zeta|,\\
K'\!\mu(z)&:=\int_\Gamma \partial_{n(z)} g(z-\zeta)\mu(\zeta)\,|d\zeta|,\quad 
W\!\mu(z):=-\partial_{n(z)}\int_\Gamma \partial_{n(\zeta)}g(z-\zeta)\mu(\zeta)\,|d\zeta|,
\end{aligned}
\end{equation}
with $g(z)={1\over2\pi}\log|z|$, then the simple layer approach implies that
upon computation of the density $\mu$, the capacity value will depend on $K'\mu$
evaluated on part of $\Gamma$, the double layer potential will involve
the estimation of $W\mu$, etc. However our approach is based on the Green formula
\begin{equation}\label{eq:bie}
V\partial_n u(z)={1\over2}u(z)+Ku(z),\quad z\in\Gamma,
\end{equation}
which operates with densities having direct interpretation of Dirichlet or Neumann data
of the harmonic function $u$. First, this means the simplest possible post-processing
for capacity computation; second, the extraction of singular terms of densities
at corner points of $\Gamma$ is also simpler.

Therefore we solve the equation \eqref{eq:bie} numerically; below we specify some
important details.
\begingroup
\begin{itemize}\itemsep0pt\parsep0pt\parskip0pt
\item We assume the contours $\Gamma_j$ be piecewise-Lyapunov \cite{vemp}.
 In the neighbourhood of a sharp corner with inner angle $\omega$, $0<\omega<2\pi$,
the function $u$ can be expanded in uniformly convergent series \cite{gris}
$$
u(r,\phi)=a_1r^\alpha\sin(\alpha\phi)+\sum_{k=2}^\infty a_k r^{k\alpha}\sin(k\alpha\phi),
\quad \alpha=\pi/\omega
$$
which allows to extract additive singular term at mesh steps belonging to a smooth part of $\Gamma$ incident to the sharp corner.
For cusps, less is known, however the representation of $\partial_n u$ on the boundary
in terms of the distance $s$ from the turning point is available \cite{mazya}, e.g. for the inner cusp one has
$$
\partial_{n_\pm} u(s)=b_0s^{-1/2}\pm b_1\log{1\over s}\pm b_2+O\left(s^{1/2}\log{1\over s}\right),
$$
which directly gives additional terms for the approximation of the solution at mesh steps near the singularity.
These expansions are used to enrich the linear space of basis functions used to represent the sought density.
\item We suppose that after subtraction of certain linear combination of
 `singular' basis functions supported on smooth parts of $\Gamma$,
the remaining part of density is periodic and smooth, so that if each $\Gamma_j$ is parametrized by scalar
$t\in[0,2\pi[$, that part can be approximated by trigonometric functions $e^{imt}$ with $|m|<N_j$ for some positive integer $N_j$.
\item The right-hand side of \eqref{eq:bie} involving $K$ acting on constant-valued Dirichlet data can be computed analytically
using the sum of angles enclosing the boundary parts with account of signs of these parts (Gauss' formula \cite{vemp}).
\item The action of $V$ on `singular' basis functions can be expressed using elementary functions.
\item The kernel of $V$ has an integrable singularity at each value of source parameter $s$ where the contour $\Gamma_j$
coincides with the observation point corresponding to parameter $t$. Were $\Gamma_j$ a circle, the respective kernel
would be (up to a multiplicative constant) equal to $\log\left|\sin(t-s)/2\right|$, the Symm kernel \cite{proess}.
Therefore for each value of $t$ we define the solutions of equation $\Gamma_j(t)=\Gamma_j(t+\Delta s)$ and
subtract from the kernel of $V$ the shifted Symm kernel $\log\left|\sin(t-s+\Delta s)/2\right|$. Thus we obtain
a smooth periodic kernel; its action on trigonometric function is approximated by rectangular quadrature rule.
\item For the integral operator with Symm's kernel, our trigonometric functions are eigenfunctions and the
eigenvalues are known \cite{proess}.
\item We use collocation at centers of mesh steps which in considered configurations can lead to linear systems
with rectangular matrices since number of mesh steps $M$ for symmetric contours is even while $\sum_jN_j$ can be odd. 
In this case, the resulting equations are understood in the least-square sense, and since the condition number 
is $O(M^2)$ \cite{schwab} it is better to use the preconditioner corresponding to described kernel splitting.
Up to matrix orders of 10000, the number of preconditioned GMRES (or CG) iterations is between 10 and 30.
\end{itemize}
\endgroup

We describe briefly the post-processing step. Recall that each $\Gamma_j$ belongs to either of two terminal
sets $\partial_0{\cal A}$,
     $\partial_1{\cal A}$,
and consider 2 Dirichlet problems:
\begin{equation}\label{eq:cap_basis}
u^{(k)}\big|_{\partial_j{\cal A}}=\delta_{kj},\quad k,j=1,2.
\end{equation}
We use the Green formula
\begin{equation}\label{eq:greencap}
\int_{\cal A} \left((\nabla u)^T\nabla w+w\Delta u\right)\,dxdy=
\int_{\partial_0{\cal A}\cup\partial_1{\cal A}}w\,\partial_n u\,|dz|
\end{equation}
for $u$, $w$ spanning the set $U=\{u^{(1)}, u^{(2)}, 1\}$. 
We see that
the Dirichlet energy of $u(z)=a_1u^{(1)}(z)+a_2u^{(2)}(z)$ is the quadratic form
$$
a^T\left[\int_{\partial_j{\cal A}}\partial_n u^{(k)}\,|dz|\right]_{k,j=1}^2a,\qquad a^T=[a_1\bis a_2].
$$
Finally, plugging $u,w\in U$ into \eqref{eq:greencap} we see that the matrix of this form
is symmetric and that the vector $q=[1\bis 1]^T$ belongs to its nullspace. Therefore this matrix is
$\lambda q_Oq_O^T$ where $q_O=[1\bis-1]^T/\sqrt2$,
%$${\lambda\over2}\begin{pmatrix}1 & -1\\ -1 & 1\end{pmatrix}
%$$
so that for $u=u^{(1)}-u^{(2)}$ the charges at the terminals are $\pm\lambda$ while the potential
difference is $2$, hence every entry of the considered matrix (by modulus) gives the capacity
of the condenser with terminals $\partial_s{\cal A}$, $s=1,2$.
Therefore capacity calculation amounts to approximating the integral of the solution
of \eqref{eq:bie} over a terminal.

\section{Numerical experiments}
\subsection{Condenser configurations}
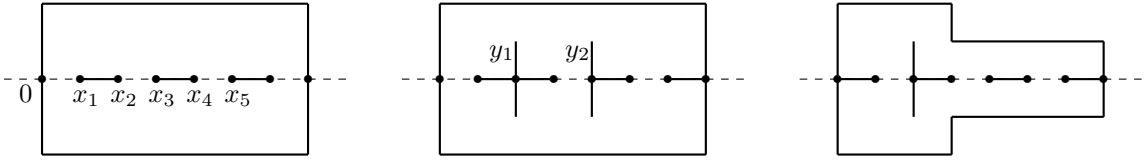
\begin{figure}[t]
\vspace*{-1cm}
\begin{picture}(40, 30)
\linethickness{.5mm}
\thinlines
\multiput(0,10)(2,0){23}{\line(1,0){1}}

\thicklines
\multiput(10,10)(10,0){3}{\line(1,0){5}}
\multiput(5,10)(5,0){8}{\circle*{1}}
\put(5,0){\line(1,0){35}}
\put(5,20){\line(1,0){35}}
\put(5,0){\line(0,1){20}}
\put(40,0){\line(0,1){20}}
\put(9,7){$x_1$}
\put(14,7){$x_2$}
\put(19,7){$x_3$}
\put(24,7){$x_4$}
\put(29,7){$x_5$}
\put(2,7){$0$}
\end{picture}
\hspace{1cm}
\begin{picture}(40, 30)
\linethickness{.5mm}
\thinlines
\multiput(0,10)(2,0){23}{\line(1,0){1}}

\thicklines
\multiput(15,10)(10,0){3}{\line(1,0){5}}
\multiput(5,10)(5,0){8}{\circle*{1}}
\put(10,10){\line(1,0){5}}
\put(5,0){\line(1,0){35}}
\put(5,20){\line(1,0){35}}
\put(5,0){\line(0,1){20}}
\put(40,0){\line(0,1){20}}
\put(15,5){\line(0,1){10}}
\put(25,5){\line(0,1){10}}
\put(11.5,13){$y_1$}
\put(21.5,13){$y_2$}

\end{picture}
\hspace{1cm}
\begin{picture}(40, 30)
\linethickness{.5mm}
\thinlines
\multiput(0,10)(2,0){23}{\line(1,0){1}}

\thicklines
\multiput(5,10)(10,0){4}{\line(1,0){5}}
\multiput(5,10)(5,0){8}{\circle*{1}}

\put(5,0){\line(1,0){15}}
\put(5,20){\line(1,0){15}}
\put(20,5){\line(1,0){20}}
\put(20,15){\line(1,0){20}}

\put(40,5){\line(0,1){10}}
\put(20,0){\line(0,1){5}}
\put(20,15){\line(0,1){5}}
\put(5,0){\line(0,1){20}}

\put(15,5){\line(0,1){10}}
\end{picture}
\caption{ Condenser configurations A, B, C (left to right).  Dashed line is the symmetry axis.}
\label{ABC}
\end{figure}

\begin{figure}[t]
\linethickness{.8mm}

\begin{picture}(35, 30)
\thinlines
\multiput(0,10)(2,0){18}{\line(1,0){1}}

\thicklines
\multiput(10,10)(10,0){2}{\line(1,0){5}}
\multiput(5,10)(5,0){6}{\circle*{1}}

\put(5,5){\line(1,0){10}}
\put(5,15){\line(1,0){10}}

\put(15,25){\line(1,0){10}}
\put(15,-5){\line(1,0){10}}

\put(25,20){\line(1,0){5}}
\put(25,0){\line(1,0){5}}

\put(30,0){\line(0,1){20}}
\put(5,5){\line(0,1){10}}
\put(15,-5){\line(0,1){10}}
\put(15,15){\line(0,1){10}}

\put(25,20){\line(0,1){5}}
\put(25,-5){\line(0,1){5}}
\put(2,12){$l_1$}
\put(9,16){$l_2$}
\put(16,19){$l_3$}

\end{picture}
\hspace{1cm}
\begin{picture}(15, 30)
\thinlines
\multiput(0,10)(2,0){10}{\line(1,0){1}}
\multiput(5,10)(10,0){2}{\circle*{1}}
\thicklines
\put(5,0){\line(0,1){20}}
\put(15,5){\line(0,1){10}}
\end{picture}\hspace{1cm}
\begin{picture}(35, 30)
\linethickness{.5mm}
\thinlines
\multiput(0,10)(2,0){20}{\line(1,0){1}}
\multiput(20,-7)(0,2){18}{\line(0,1){1}}
\thicklines
\put(5,-5){\line(1,0){30}}
\put(5,25){\line(1,0){30}}
\put(10,5){\line(1,0){20}}
\put(10,15){\line(1,0){20}}

\put(5,-5){\line(0,1){30}}
\put(35,-5){\line(0,1){30}}
\put(10,5){\line(0,1){10}}
\put(30,5){\line(0,1){10}}
\end{picture}
\hspace{1cm}
\begin{picture}(35, 30)
\thinlines
\multiput(0,10)(2,0){20}{\line(1,0){1}}
\multiput(20,-7)(0,2){18}{\line(0,1){1}}

\thicklines
\put(5,-5){\line(1,0){30}}
\put(5,25){\line(1,0){30}}
\put(5,-5){\line(0,1){30}}
\put(35,-5){\line(0,1){30}}

\put(10,5){\line(1,0){5}}
\put(10,15){\line(1,0){5}}
\put(10,5){\line(0,1){10}}
\put(15,15){\line(0,1){5}}
\put(15,0){\line(0,1){5}}
\put(15,20){\line(1,0){10}}
\put(15,0){\line(1,0){10}}
\put(25,15){\line(0,1){5}}
\put(25,0){\line(0,1){5}}
\put(25,5){\line(1,0){5}}
\put(25,15){\line(1,0){5}}
\put(30,5){\line(0,1){10}}
\end{picture}
\vspace{4mm}
\caption{ Condenser configurations D, E, F, G (left to right).  Dashed line is the symmetry axis.}
\label{DEFG}
\end{figure}
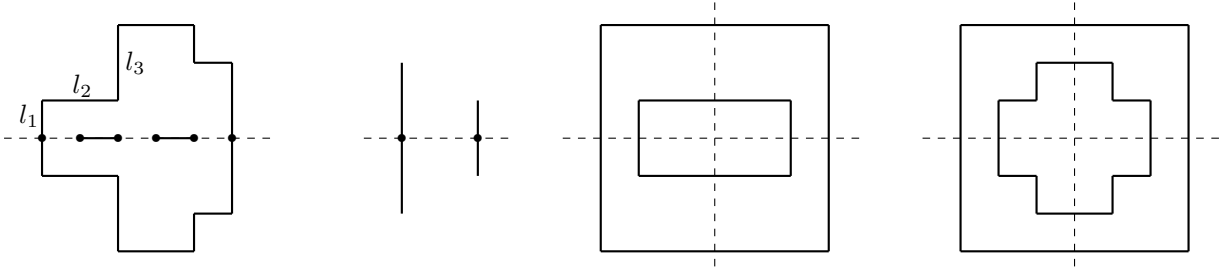

In this section we describe seven topological types  -- or configurations -- of axisymmetric condensers. We designate them 
as A,B,\dots, G, as shown in Figs.~\ref{ABC}, \ref{DEFG}.  

\subsection{Dimensions of condensers and their capacities}
We need a unified system for encoding the  dimensions of the multiply connected polygons we mentioned above.
For configurations A,B,C,D,E (see Figs.~\ref{ABC}, \ref{DEFG})
we first describe all significant points on the (symmetry) real axis: the
endpoints of the horizontal slots, the bases of vertical slots and also
the intersection points of the outer contour and the symmetry axis. Those points are distinguished in the pictures. 
We put  all the numbers in increasing order to the array X, starting with zero (which means picture
normalization). For instance, it is possible that  $X=[0, 1, 2, 3, 4, 5, 6, 7]$ for configurations A,B,C of Fig.~\ref{ABC}.
For the configurations D,E of Fig.~\ref{DEFG} the lengths of array $X$ is 6 and 2, respectively.

Next, we  have to codify the heights of the vertical cuts emanating
from the real axis (left to right). We put them to the array $Y$, for instance
in cases A and D the array $Y$ is empty (no vertical cuts), and in cases B,C it has the length 2 and 1, respectively.

And finally we need to describe the shape of the outer contour.This can be done
by specifying the lengths of the consecutive segments of the chain in
the upper half plane (with the exception of the last two which give no
new information),  starting from the origin which we put to the left
intersection point of the outer contour with the real axis. We put
the lengths to the array $L$. For instance, $\sharp L=2$ in cases A,B, $\sharp L=3$  and 4 
for configurations C,D, respectively and this array is empty in case E. 

For the remaining shapes F,G of Fig.~\ref{DEFG}: the concentric rectangles and  the
`swiss cross' we describe the the inner and the outer contour in a similar manner.  
First  separate  the  North-West
quarter of the picture  and encode the arising chains as previously.
Now we have two arrays $L$ instead of one: $L1$ for the outer contour
and  $L2$ for the inner one. For instance, both arrays contain two elements for the shape F, while  
for the shape  G the first array $L1$ contains 2 elements and the second  -- 4.

Having adopted this system of designation,  we formulate the numerical experiments and obtain the capacities  
in the form of the tables below.
\begin{table}[h!] 
{\small
\begin{tabular}{c|l|l|l}
Id.& X & L& Capacity\\
\hline
A1&0,1,3,4,5,6,9,11 &2 & 9.72079120617096926\\
A2&0,2,2.5,4.5,5,6,10.5,11  &1 & 15.8964033734093744\\
A3&0,.5,5,7.5,8,10,10.5,11 &1 & 16.5708349921371510\\
A4&0,1,1.5,2,8,10,10.5,11 & 3& 9.0776774351927967\\
A5&0,2,2.5,7.5,8,9,10.9,11  &1 & 12.1642765126444534 \\
A6& 0,3,5,6,8,9,10,11 &5 & 5.59517889911177450\\
\end{tabular}
\hfill
\begin{tabular}{c|l|l|l|l}
Id.& X & Y& L& Capacity\\
\hline
B1&0,1,2,3,5,6,7,8&2,1&3&8.86185570899657537\\
B2&0,1,2,3,4,5,6,7&1,1&2&8.28583441065142426\\
B3&0,1,2,3,4,5,6,7&2,1&3&8.22274382175325185 \\
B4&0,1,2,3,4,5,6,7&1,2&3&8.11029353036022815\\
B5&0,1,2,6,7,9,10,12&1,1&2&12.17857832040164176 \\
B6&0,1,2,6,7,9,10,11&1,3&4&10.31120091451990165\\
\end{tabular}}
\caption{Dimensions for configuration A (left) and B (right) condenser and their capacities}
\end{table}

\begin{table}[t!]
\begin{tabular}{c|l|l|l|l}
Id.&Y&X&L& Capacity\\
\hline
C1&2&0,1,2,3,5,6,7,8&3,3,1&9.438272363758330697\\
C2&1&0,1,2,3,4,6,7,8&2,3,1&11.027047279861000458\\
C3&2&0,1,2,3,4,5,7,8&3,3,1&8.777515831134811065\\
C4&3&0,1,2,3,4,5,6,8&4,3,3&14.988032697667965659\\
C5&3&0,2,3,4,5,6,7,8&4,5,3&11.391736530234725936\\
C6&3&0,1,2,4,5,6,7,8&4,7,2&9.282399749809620850\\
\end{tabular}
\caption{Dimensions for configuration C condenser with corresponding capacities}
\end{table}

\begin{table}[t]
\begin{tabular}{c|l|l|l}
Id.& X & L& Capacity\\
\hline
D1&0,1,3,4,5,8&2,3,2,2,3&6.298067056123278293\\
D2&0,2,3,4,7,8&2,3,2,2,3&8.994834022659064427\\
D3&0,4,5,6,7,8&2,3,2,2,3&6.450372178406949499\\
D4&0,1,3,5,7,8&2,3,3,2,3&7.438309246517998359\\
D5&0,1,3,4,5,8&2,5,1,2,1&5.801923089413399328\\
D6&0,1,3,6,7,8&2,3,4,2,5&7.753127034648571466\\
\end{tabular}
\hfill
\begin{tabular}{c|l|l}
Id.& Y & Capacity\\
\hline
E1&1,2&1.56994325474948999\\
E2&2,2&1.87306699654806386\\
E3&3,2&2.08203777712328096\\
E4&4,2&2.23259828277206300\\
E5&5,2&2.34158897620030515\\
E6&3,3&2.35241226225174034\\
\end{tabular}
\caption{Dimensions for configuration D condenser (left)
and  E  condenser `two vertical slots` (right) with $X=[0,5]$ and their capacities}
\end{table}

\begin{table}[t]
\begin{tabular}{c|l|l|l}
Id.& L1 & L2& Capacity\\
\hline
F1&3,4&1,1&5.6327570222823258486\\
F2&3,4&.3,3&8.5383099064779181521\\
F3&3,4&2,.1&5.7845537023572573861\\
F4&1,4&.5,3&28.5499438953187884\\
F5&2,5&1,4&22.234504016933507380\\
F6&2,7&1.5,.2&7.6417584869737709307\\
\end{tabular}
\hfill
\begin{tabular}{c|l|l|l}
Id.& L1 & L2& Capacity\\
\hline
G1&4, 5&1,2,1,1&9.578338769355109451\\
G2&4, 4&1,2,2,1&16.076240045355723868\\
G3&5, 4&1,2,2,1&13.192681030463681933\\
G4&4, 6&1,3,2,1&14.350501722644794417\\
G5&4, 5&1,3,2,1&17.116438880060405780\\
G6&4, 6&1,3,2,2&21.116597096285347718\\
\end{tabular}
\caption{Dimensions for configuration F condenser `concentric rectangles` (left) and G: `swiss cross` (right) and their capacities}
\end{table}

{\bf Remarks} 1) Computations for the above tables were made with double and sometimes quadruple precision and give consistent results for 
all four methods within the accuracy we guarantee. 2) Our implementation of the boundary integral equation method gave 
ca. 4-6 true digits for the computations with double accuracy. 3) The hp-FEM gave 13-14 true digits.
4) Configuration  $E$ `two slots' was not computed with $hp$-FEM since it is unbounded.
5) The theta functions and  BEM algorithms were implemented in FORTRAN, $hp$-FEM used Mathematica, Lauricella functions 
were computed with C++. 

\section{Conclusion}
Condensers occur in several problems of potential theory, see Polya-Szeg\H o \cite{posz}, and geometric function theory and thus many authors have studied capacities either theoretically or from pragmatic point of view.
Capacity computation belongs to the area of computational potential theory.
Engineering applications such as electronics, integrated circuit design often
use application specific numerical approximation methods, see e.g.

{\tt http://ieeexplore.ieee.org/abstract/document/817589/},  

\noindent
or books \cite{posz,sl}. While several such approximation methods exist, it is rare to find
systematically tabulated and double-checked numerical values of specific
condenser capacities in the current literature. Thus, taking into account the needs
of various application areas and the currently available computer facilities,
there is a lot of room for further research in computational potential theory.

The uncertainty principle applied to the numerical methods divides the latter very roughly 
into two groups: `more flexible' and `more accurate'.  The first group has more real life applications and the second may be used 
to validate, verify, tune in  and calibrate first group methods.  The goal of this paper is to compute the capacities for a wide set of 42 axially symmetric planar condensers by two independent high precision analytical computational methods. One of those  uses Lauricella function, the generalization of hypergeometric  one; the other employs (multidimensional) theta functions introduced by Jacobi, Rosenhein, Goeppel and Riemann. Having done this we analyze the performance of more traditional computational approaches  -- \emph{hp-Finite Element Method} \cite{hrv1, hrv2, hrv3}  and \emph{Boundary Integral Equations Method}. All these methods have been previously
reported and tested in the literature.

In this paper we have used the simplest concept of capacity which comes from 
physics or electrical engineering as a factor relating the voltage of the condenser to the electrical charge on its plates. 
This concept delivers a variety of mathematically significant problems. Other settings also deserve further investigation:

{\bf Capacity matrix.}~
Given a multiply connected domain $\cal A$ in the Riemann sphere with a `good' boundary, one can apply locally constant voltage  to all its boundary components. The arising electric potential is merely a harmonic function $u(z)$ with locally constant boundary values. The Dirichlet integral of this function calculated over the domain becomes a quadratic form of the values of the potential at the boundary components. The matrix of this quadratic form is called the capacity matrix. It has the function theoretic meaning
\cite{Bt} as the period matrix of M-curve, the double of the planar domain $\cal A$. In particular,
it is symmetric and diagonally dominant.   
 
{\bf Logarithmic capacity.}~
In case the domain $\cal A$ contains the infinity inside, one can consider its Green's function.
Namely, the harmonic function with logarithmic singularity at infinity and vanishing at the boundary of $\cal A$. The constant term in the decomposition of the Green's function gives the log capacity $C({\cal A})$ of $\cal A$:  $G(z)=\log(|z|)-\log(C({\cal A}))+o(1)$, $z\to\infty$.

\section*\refname


\small
\begin{thebibliography}{99}

\bibitem[ABR]{abr}
\textsc{S. Axler, P. Bourdon, and W. Ramey,} Harmonic Function Theory. 
2nd ed., Springer, 2001.

\bibitem[AB]{ab}
\textsc{L.V. Ahlfors and A. Beurling,} Conformal invariants 
and function theoretic null sets. Acta Math. 83 (1950), 101--129.

\bibitem[AC]{ac}
\textsc{J. Albrecht and L. Collatz, eds.,} 
Numerical Treatment of Integral Equations, Birkh\"auser, Basel, 1980.

\bibitem[AVV]{avv}
\textsc{G. D. Anderson, M. K. Vamanamurthy, and M. Vuorinen,} 
Conformal invariants, inequalities and quasiconformal maps. Wiley, 1997.

\bibitem[Ah]{ah} 
\textsc{L. V. Ahlfors,} Conformal invariants: Topics in
Geometric Function Theory. McGraw-Hill, New York, 1973.

\bibitem[A]{A} 
\textsc{N. I.  Akhiezer,} 
Elements of the theory of elliptic functions. 
AMS, Providence, RI, 1990.


\bibitem[B]{B} 
\textsc{A. B. Bogatyr\"ev,} The conformal mapping of rectangular heptagons.
Sb. Math. 203 (2012), no. 11-12, 1715--1735.

\bibitem[B2]{B2} 
\textsc{A. B. Bogatyr\"ev,}
Image of Abel-Jacobi map for hyperelliptic genus 3 and 4 curves // Journal Approx. Theory, 191 (2015), 38--45.

%\bibitem[BE]{be}
%\textsc{W. Bergweiler and A. Eremenko,} Goldberg's constants. 
%J. Anal. Math. 119 (2013), no. 1, 365--402.

\bibitem[BG1]{BG1} 
\textsc{A. B. Bogatyr\"ev and O. A. Grigor'ev,}  Closed formula for the capacity of several aligned segments, 
Proceedings of Steklov Inst.,  298 (2017), 67--74.

\bibitem[BG2]{BG2} 
\textsc{A. B. Bogatyr\"ev and O. A. Grigor'ev,}
Conformal mapping of rectangular heptagons II,
Published online in Comp. Meth. Func. Theory, DOI:10.1007/s40315-017-0217-z,
arXiv:1612.01127.

\bibitem[Bt]{Bt} 
\textsc{A. B. Bogatyrev,}
Real meromorphic differentials: A language for describing meron configurations 
in planar magnetic nanoelements // Theor. and Math. Phys., 193:1 (2017), 1547–1559; arXiv: 1610.04984.


\bibitem[BSV]{bsv} 
\textsc{D. Betsakos, K. Samuelsson, and M. Vuorinen,} The computation of capacity 
of planar condensers. Publ. Inst. Math. (Beograd) (N.S.) 75(89) (2004), 233--252.

%\bibitem[BV]{bv} 
%\textsc{D. Betsakos and M. Vuorinen, } Estimates for conformal capacity. Constr. 
%Approx. 16 (2000), no. 4, 589--602.


\bibitem[Be1]{be1} 
\textsc{S. I. Bezrodnykh,}
Analytic continuation formulas and Jacobi--type relations for Lauricella function.
Dokl. Math. 93 (2016), no. 2, 129--134.

\bibitem[Be2]{be2} 
\textsc{S. I. Bezrodnykh,} On the analytic continuation of the Lauricella function 
$F_D^{(N)}$. Math. Notes. 100 (2016), no. 1--2, 318--324.

\bibitem[Be3]{be3} 
\textsc{S.I. Bezrodnykh,} 
Lauricella hypergeometric function $F_D^{(N)}$, the Riemann-Hilbert problem and some applications// Russian Mathematical Surveys. 73 (2018), no. 6, 941–1031.

\bibitem[BeV1]{bevl1} 
\textsc{S.I. Bezrodnykh, V.I. Vlasov,}
The Riemann--Hilbert problem in a complicated domain for the model of magnetic reconnection
in plasma. Comput. Math. Math. Phys. 42 (2002), no 3, 263--298.

\bibitem[BeV2]{bevl2} 
\textsc{S.I. Bezrodnykh, V.I. Vlasov,}
Singular Riemann--Hilbert problem in complex-shaped domains.
Comput. Math. Math. Phys. 54 (2014), no 12, 1826--1875.

\bibitem[Bo]{bo}
\textsc{F. Bosshard,} Die Konstruktion konformer Abbildungen mit der Methode
der finiten Elemente, Dissertation, Z\"urich, 1980.


\bibitem[CWLH]{cwlh}
\textsc{S. N. Chandler-Wilde, J. Levesley, and D. M. Hough,} 
Evaluation of a boundary integral representation for the conformal mapping of 
the unit disk onto a simply-connected domain, Adv. Compt. Math. 3 (1995), 115--135.

\bibitem[Cr]{cr}
\textsc{D. Crowdy,} The Schwarz--Christoffel mapping to bounded multiply 
connected polygonal domains, Proc. R. Soc. Lond. Ser. A Math. Phys. Eng.Sci, 461. 
(2005) 2653--2678.

\bibitem[DrT]{dt} 
\textsc{T. A. Driscoll and L. N. Trefethen,} Schwarz--Christoffel mapping. 
Cambridge University Press, 2002.

\bibitem[Da]{da}
\textsc{H. Daeppen, } Die Schwarz-Christoffel-Abbildung f\"ur zweifach zusammenh\"angende 
Gebiete mit Anwendungen. Ph.D. thesis, E.T.H., Z\"urich, 1988.

\bibitem[DeK]{dek}
\textsc{T.K. DeLillo, E.H. Kropf,} Numerical computation of the Schwarz--Christoffel 
transformation for multiply connected domains, SIAM J. Sci. Comput. 33 (2011) 1369--1394.

\bibitem[DeEK]{deek}
{\sc T. DeLillo, A. Elcrat, and E. Kropf,} Calculation of resistances for multiply connected domains using Schwarz-Christoffel
transformations, Comput. Methods Funct. Theory, 11 (2011), 725-745.


\bibitem[Dr]{dr}
\textsc{T.A. Driscoll,} Schwarz--Christoffel toolbox for MATLAB, 
\url{http://www.math.udel. edu/~driscoll/SC/}

\bibitem[Du]{du} 
\textsc{V. N. Dubinin,} Condenser Capacities and Symmetrization in Geometric 
Function Theory, Birkh\"auser, 2014.

\bibitem[E]{e} 
\textsc{H. Exton,}
Multiple hypergeometric functions and application.
J. Wiley \& Sons, New York, 1976.

\bibitem[ET]{et}
{\sc M. Embree and L. N. Trefethen,} Green's functions for multiply connected domains via
conformal mapping, SIAM Review, 41 (1999),  721-744.

\bibitem[Ga64]{ga64} 
\textsc{D. Gaier,} Konstruktive Methoden der konformen Abbildung. 
Springer-Verlag, 1964.

\bibitem[Ga72]{ga72}
\textsc{D. Gaier,} 
Ermittlung des konformen Moduls von Vierecken mit Differenzenmethoden. 
Numer. Math. 19 (1972) 179-194.

\bibitem[Ga95]{ga95} 
\textsc{D. Gaier,} Conformal modules and their computation. In: Computational 
Methods and Functional Theory (CMFT'94). Eds.: R.M.Ali et al. World Scientific, 1995. 
159--171. 

\bibitem[GaH]{gah} 
\textsc{D. Gaier, W.K. Hayman,} 
On the computation of modules of long quadrilaterals. Constr. Approx. 7 (1991), 459--467.

\bibitem[Gd]{gris} 
\textsc{P. Grisvard,}
Elliptic Problems in Nonsmooth Domains.
Pitman Publ. Ltd., London, 1985. 410pp.

\bibitem[Gr]{gr} 
\textsc{ O. A. Grigoriev,} Numerical-analytical method for conformal mapping 
of polygons with six right angles, Comp. Math. Math. Phys., 53:10, 1629--1638 (2013)

\bibitem[HKW]{hkw}
\textsc{G. C. Hsiao, P. Kopp, and W. L. Wendland,} A Galerkin collocation method 
for some integral equations of the first kind, Computing 25 (1980), 89--130.

\bibitem[HLCW]{hlcw}
\textsc{D. M. Hough, J. Levesley, and S. N. Chandler--Wilde,} 
Numerical conformal mapping via Chebyshev weighted solutions of Symm's integral equation, 
J. Comput. Appl. Math. 46 (1993), 29--48.

\bibitem[HPS]{hps}
\textsc{J. Hersch, A. Pfluger, A. Schopf,} \"Uber ein
    simultanes Differenzenverfahren zur Absch\"atzung der
    Torsionssteifigkeit und der Kapazit\"at nach beiden Seiten, 
    Z. Angew. Math. Physik. 7 (1956), 89--113.

\bibitem[HP]{hp}
\textsc{D. M. Hough and N. Papamichael,} The use of splines and singular functions 
in an integral equation method for conformal mapping, Numer. Math. 37 (1981), 133--147.

\bibitem[HRV1]{hrv1} 
\textsc{H. Hakula, A. Rasila and M. Vuorinen,} On moduli of rings and quadrilaterals: 
algorithms and experiments. SIAM Sci. Comput. 33 (2011), no. 1, 279--302. 

\bibitem[HRV2]{hrv2} 
\textsc{H. Hakula, A. Rasila and M. Vuorinen,} Computation of exterior moduli 
of quadrilaterals. Electron. Trans. Numer. Anal. 40 (2013), 436--451.


\bibitem[HRV3]{hrv3} \textsc{ H. Hakula, A. Rasila, and M. Vuorinen},
{Conformal modulus on domains with strong singularities and cusps}.- { Electron. Trans. Numer. Anal. 48, 462--478, 2018,} DOI: 10.1553/etna\_vol48s462, { arXiv:1501.06765[math.NA]}

\bibitem[H]{h}
\textsc{D. M. Hough,} Jacobi polynomial solutions of first kind integral
equations for numerical conformal mapping. J. Comput. Appl. Math. 13 (1985) 359--369.

\bibitem[HaM]{ham}
\textsc{H. Hara and H. Mizumoto}, Determination of the modulus of
quadrilaterals by finite element methods. J. Math. Soc. Japan. 42 (1990), 295--326.

\bibitem[Ha17]{Hakula2017} 
\textsc{H. Hakula, M. Neilan, and J. S. Ovall},
A Posteriori Estimates Using Auxiliary Subspace Techniques, J. Sci. Comp. 72 (1), 2017, 97--127.

\bibitem[He]{he} 
\textsc{P. Henrici,} Applied and computational complex analysis.
Vol.~3. Wiley Classics Library, 1993.

\bibitem[HoT]{hot}
\textsc{L. H. Howell and L.N. Trefethen,} A modified
    Schwarz--Christoffel transformation for elongated regions. SIAM J. 
    Sci. Stat. Comp. 11 (1990), 928--949.

\bibitem[Hu]{hu}
\textsc{C. Hu,} A software package for computing Schwarz--Christoffel
conformal transformation for doubly connected regions. ACM Trans. Math. 
Soft. 24 (1998), 317--333.

\bibitem[Hyv]{hyv}
\textsc{N. Hyv\"onen,} Complete electrode model of electrical impedance tomography:
approximation properties and characterization of inclusions. SIAM J. Appl. Math. 64
(2004), 902--931.

\bibitem[JS]{js} 
\textsc{M.A. Jaswon, G.T. Symm,} 
Integral equation methods in potential theory and elastostatics. London:
Academic Press, 1977.

\bibitem[Kir]{kir} 
\textsc{Kirsch, Siegfried}, Transfinite diameter, Chebyshev constant and capacity. Handbook of complex analysis: geometric function theory. Vol. 2, 243–308, Elsevier Sci. B. V., Amsterdam, 2005. 

\bibitem[KS]{ks} 
\textsc{ W. von Koppenfels and F. Stallman}, 
Praxis der konformen Abbildung. Springer-Verlag, 1959.

%\bibitem[Klo]{klo} 
%\textsc{H. Kloke,} On the capacity of a plane condenser and
%conformal mapping. J. Reine Angew. Math. 358 (1985), 179--201.

\bibitem[KrSg]{krut} 
\textsc{P. A. Krutitskii, A. I. Sgibnev,}
The method of integral equations in the generalized jump problem for the Laplace equation
outside cuts on the plane. Differ. Equ. 38 (2002), no. 9, 1277--1292.

\bibitem[Ku]{ku}
\textsc{ R. K\"uhnau,} (Ed.) Handbook of Complex Analysis: Geometric Function Theory. 
Vol. 2. Amsterdam: North Holland/Elsevier, 2005.

\bibitem[Ky]{ky}
\textsc{ P.K. Kythe,} Computational conformal mapping. Birkh\"auser, 1998.

\bibitem[LHCW]{lhcw}
\textsc{J. Levesley, D. M. Hough, and S. N. Chandler-Wilde,} A Chebyshev 
collocation method for solving Symm's integral equation for conformal mapping: a partial 
error analysis, IMA J. Numer. Anal. 14 (1994), 57--79.

\bibitem[LL]{ll}
\textsc{P.A.A. Laura,   L.E. Luisoni,} An application of conformal mapping to 
the determination of the characteristic impedance of a class of coaxial systems. 
IEEE Trans MTT. 25 (1977), 162--164.

\bibitem[LPS]{lps} 
\textsc{D. Levin, N. Papamichael, A. Sideridis,} 
The Bergman kernel method for numerical conformal mapping of simply connected
domains. J. Inst. Math. Appl. 22 (1978), 171--187.

\bibitem[LSN]{lsn} 
\textsc{J. Liesen, O. S\`ete and  M. M.S. Nasser,}
Fast and accurate computation of the logarithmic capacity of compact sets,
Comp.  Methods and Function Theory
Dec. 2017, Vol17, Issue 4, 689--713.
arXiv:1507.05793[math.NA]

\bibitem[L]{l}
\textsc{R.I. Laugesen,} Extremal problems involving logarithmic and Green 
capacity. Duke Math. J. 70 (1993), 445--480.

\bibitem[MSS]{mss} 
\textsc{V. V. Maymeskul, E. B. Saff, and N. S. Stylianopulos,} 
$L^2$-Approximations of power and logarithmic functions with applications to 
numerical conformal mapping. Numer. Math. 91 (2002), 503--542.

%\bibitem[Men]{men} 
%\textsc{K. Menke,} Point systems with extremal properties and conformal
%mapping. Numer. Math. 54 (1988), 125--143.

\bibitem[Mi]{mi}
\textsc{H. Mizumoto,} An application of Green's formula of a discrete function: 
Determination of periodicity moduli I, II. Kodai Math. Sem. Rep. 22 (1970), 244--249.

\bibitem[M]{Mum}  
\textsc{D. Mumford,} Tata lectures on Theta. Vol. I--II, Springer, 1983.

\bibitem[MaS]{mazya} 
\textsc{V. G. Mazya, A. A. Soloviev,}
On an integral equation for the Dirichlet problem in a plane domain with cusps on the boundary. 
Math. USSR Sb. 68 (1991), no. 1, 61--83.

\bibitem[O67]{o67} 
\textsc{G. Opfer,} Untere, beliebig verbesserbare Schranken f\"ur den
Modul eines zweifach zusammenh\"angenden Gebietes mit Hilfe von
Differenzenverfahren, Dissertation, Hamburg, 1967.

\bibitem[P03]{p03}
\textsc{N. Papamichael,} Dieter Gaier's contribution to conformal mapping. 
Comp. Meth. Funct. Theory. 3 (2003) 1--53.

\bibitem[P89]{p89} 
\textsc{N. Papamichael,} Numerical conformal mapping onto a rectangle with 
application to the solution of Laplace problems. J. Comput. Appl. Math. 28 (1989),
63--83.

%\bibitem[PK]{pk}
%\textsc{N. Papamichael and C.A. Kokkinos,} The use of singular
%   functions for the approximate conformal mapping of doubly-connected
%   domains. SIAM J. Sci. Stat. Comp. 5 (1984) 684--700.

\bibitem[PPSS]{ppss}
\textsc{N. Papamichael, I. E. Pritsker, E. B. Saff, and N. Stylianopoulos,} 
Approximation of conformal mappings of annular regions. Numer. Math. 76 (1997), 
489--513.

%\bibitem[PRS]{prs}
%\textsc{N. Papamichael, S. Ruscheweyh, and E. B. Saff, eds.,} 
%Computation Methods and Function Theory, (CMFT'97), World Scientific, 1999.

\bibitem[PS]{ps} 
\textsc{N. Papamichael and N. Stylianopoulos,} Numerical conformal mapping. 
Domain decomposition and the mapping of quadrilaterals. World Scientific Publishing 
Co. Pte. Ltd., Hackensack, NJ, 2010. xii+229 pp.

\bibitem[PWH]{pwh} 
\textsc{N. Papamichael, M. K. Warby, and D.M. Hough}, The treatment of corner
and pole-type singularities in numerical conformal mapping technique.
    J. Comp. Appl. Math. 14 (1986), 163--191.

\bibitem[PtSw]{schwab} 
\textsc{T. von Petersdorff, C. Schwab,}
Wavelet approximations for the first kind boundary integral equations on polygons.
Numer. Math. 74 (1996), no. 4, 479--516.

\bibitem[PrSi]{proess} 
\textsc{S. Pr\"ossdorf, B. Silbermann,}
Numerical analysis for integral and related operator equations.
Birkh\"auser, Basel, 1991.

\bibitem[PoSz]{posz} 
\textsc{G. Polya and G. Szeg\"o,} Isoperimetric inequalities in mathematical 
physics. Princeton,  Univ. Press, 1951.

\bibitem[RC]{rc}
\textsc{ G.J. Rogers, G.K. Cambell,} The piece-by-piece solution of elliptic 
boundary value problems. J. Phys. D: Appl. Phys. 8 (1975), 1615--1623.

\bibitem[RRR]{rrr} 
\textsc{Q. Rajon, Th. Ransford and J. Rostand,} 
Computation of capacity via quadratic programming. J. Math. Pures Appl. (9) 94 (2010), no. 4, 398--413.

\bibitem[RSSV]{rssv}
\textsc{St. Ruscheweyh, E. B. Saff, L. C. Salinas, and R. S. Varga (eds.),} 
Computational Methods and Function Theory, Lecture Notes in Mathematics, 
Springer-Verlag, Berlin, 1990.

%\bibitem[RV1]{rv1} 
%\textsc{A. Rasila and M. Vuorinen,}
%Experiments with moduli of quadrilaterals. Rev. Roumain Math. Pures Appl. 51 (2006), 747--757.
%arXiv math.NA/0703149, MR2320923.

%\bibitem[RV2]{rv2} 
%\textsc{A. Rasila and M. Vuorinen,}
%Experiments with moduli of quadrilaterals II. Proc. Internat. Conf. on Geometric Function Theory,
%Special Functions and Applications (ICGFT), J. Analysis. 15 (2007), 229--237.

\bibitem[Rei]{rei} 
\textsc{L. Reichel,} A fast method for solving certain integral equations 
of the 1st kind with application to conformal mapping, J. Comput. Appl. Math. 
14 (1986) no. 1--2, 125--142.

\bibitem[SK]{sk}
\textsc{ A. Seidl, H. Klose,} Numerical conformal mapping of a towel--shaped
region onto a rectangle. SIAM J. Sci. Stat. Comput. 6 (1985), 833--842.

\bibitem[SL]{sl} 
\textsc{ R. Schinzinger and P.A.A. Laura} Conformal mapping. Methods and 
applications.
Dover Publ., 2003.

\bibitem[SM]{sm}
\textsc{ E. Sharon, D. Mumford,} 2D--Shape analysis using conformal mapping.
Intern. J. Computer Vision. 70 (2006), No.~1, DOI:10.1007/s11263-006-6121-z.

\bibitem[Sw]{sw} 
\textsc{Ch. Schwab},
$p$- and $hp$-{F}inite {E}lement {M}ethods, Oxford University Press, 1998.

\bibitem[SzBa]{szba} 
\textsc{B. Szabo  and I. Babu\v ska},
Finite {E}lement {A}nalysis, Wiley, 1991.

\bibitem[TG]{tg} 
\textsc{T. Theodorsen, I.E. Garrick,} General potential theory of arbitrary wing
sections. NACA Rep. 452 (1933).

\bibitem[THG]{thg}
\textsc{A. Tiwarya, C. Hua and S. Ghosh,} Numerical conformal mapping method based on Voronoi
cell finite element model for analyzing microstructures with irregular heterogeneities,
Fin. Elem. Anal. Design 43 (6--7) (2007), 504--520.

\bibitem[TiC]{tic}
\textsc{ J.C. Tippet,  D.C. Chang,} Radiation characteristics of electrically 
small devices in a FEM transmission cell. IEEE Transactions on Electromagnetic Compatibility. 
18 (1976),  No.~4, 134--140.

\bibitem[Tr80]{tr80} \textsc{L.N.Trefethen}, Numerical computation of the Schwarz–Christoffel transformation, SIAM J. Sci. Stat. Comput., 1 (1980), 82–102.

\bibitem[Tr86]{tr86}
\textsc{ L. N. Trefethen,} Ed., Numerical conformal mapping. North--Holland,
Amsterdam, 1986.

\bibitem[Tr]{tr}
\textsc{ L. N.Trefethen, } Analysis and design of polygonal resistors by conformal
mapping. J. Appl. Math. Phys. (ZAMP). 35 (1984) 692--704.

\bibitem[Ts]{ts}
\textsc{ M. Tsuji,} Potential Theory in Modern Function Theory, Tokyo, 1959.

\bibitem[V]{vemp} 
\textsc{V. S. Vladimirov,}
Equations of Mathematical Physics.
Marcel Dekker, 1971. 464pp.

\bibitem[W79]{w79}
\textsc{J. Weisel,} L\"osung singul\"arer Variationsprobleme durch die Verfahren 
von Ritz und Galerkin mit finiten Elementen. Anwendungen in der konformen Abbildung.
Mitt. Math. Sem. Giessen. 138 (1979), 1--150.

\bibitem[Y]{y}
\textsc{ Y. Yan,} The collocation method for first-kind boundary integral 
equations on polygonal regions, Math. Comp. 54 (1990), 139--154.

\end{thebibliography}
\end{document}